\documentclass[12pt]{amsart}
\usepackage{amsmath,amsthm}
\usepackage{amssymb}
\usepackage{euscript}           

\newtheorem{thm}{Theorem}[section]
\newtheorem{prop}[thm]{Proposition}
\newtheorem{cor}[thm]{Corollary}
\newtheorem{lemma}[thm]{Lemma}

\theoremstyle{definition}
\newtheorem{defn}[thm]{Definition}

\theoremstyle{remark}
\newtheorem{remark}[thm]{Remark}

\numberwithin{equation}{section}

  \DeclareMathOperator{\card}{card} 
   
  \DeclareMathOperator{\coo}{coor} 
 \DeclareMathOperator{\de}{def}  
\DeclareMathOperator{\dist}{dist}   
   
  \DeclareMathOperator{\id}{id} \DeclareMathOperator{\im}{im}
   
    \DeclareMathOperator{\pr}{pr}
    
   \DeclareMathOperator{\spec}{spec}

\newcommand{\asy}{\textsf{asymp}}

\newcommand{\et}{{\textsl{\'et}}}

\newcommand{\fil}{\mathsf{fill}}

\newcommand{\isop}{\mathsf{isop}}

\newcommand{\Los}{\L\'os}

\renewcommand{\th}{{\text{th}}}

\newcommand{\topo}{{\textsl{top}}}

\newcommand{\topos}{\EuScript}
\newcommand{\cate}{\mathcal}

\newcommand{\CC}{\mathbb{C}}

\newcommand{\FF}{\mathbb{F}}
\newcommand{\G}{\topos G}

\newcommand{\tG}{\mathsf G}

\newcommand{\K}{\cate K}

\newcommand{\N}{\mathbb{N}}

\newcommand{\tP}{\mathsf P}

\newcommand{\QQ}{\mathbb{Q}}

\newcommand{\RR}{\mathbb{R}}

\newcommand{\tT}{\mathsf T}
\newcommand{\U}{\topos U}

\newcommand{\V}{\cate V}

\newcommand{\Z}{\mathbb{Z}}
\newcommand{\cZ}{\cate Z}
\newcommand{\tW}{\mathsf W}

\newcommand{\nil}{\emptyset}

\newcommand{\te}{\bd{\textup{1}}}

\newcommand{\lc}{\langle}                       
\newcommand{\rc}{\rangle}                       

\renewcommand{\leq}{\leqslant}          
\renewcommand{\geq}{\geqslant}          

\newcommand{\iso}{\approx}          
\newcommand{\tens}{\otimes}         

\newcommand{\dfra}{\Rightarrow}         
\newcommand{\dlra}{\Leftrightarrow}     

\newcommand{\epi}{\twoheadrightarrow}
\newcommand{\ra}{\rightarrow}
\newcommand{\lra}{\longrightarrow}

\newcommand{\inc}{\hookrightarrow}

\newcommand{\bd}{\textbf}

\newcommand{\llra}[1]{\stackrel{#1}{\lra}}      
\newcommand{\xra}[1]{\xrightarrow{#1}}          

\usepackage[all]{xy}

\xymatrixcolsep{1.9pc}                     
\xymatrixrowsep{1.9pc}
\newdir{ >}{{}*!/-5pt/\dir{>}}             
\newdir{ |}{{}*!/-5pt/\dir{|}}             

\begin{document}
\title{Isoperimetric inequalities and the Friedlander--Milnor conjecture}
\author{Tibor Beke}
\address{Department of Mathematics\\
University of Michigan\\
525 East University Avenue\\
Ann Arbor, MI 48109}
\email{tbeke@umich.edu}
\date{\today}
\begin{abstract}
We prove that Friedlander's generalized isomorphism conjecture on the cohomology of algebraic groups, and hence the Isomorphism Conjecture for the cohomology
of the complex algebraic Lie group $\tG(\CC)$ made discrete, are equivalent to the existence of an isoperimetric inequality in the homological bar complex of
$\tG(F)$, where $F$ is the algebraic closure of a finite field.
\end{abstract}
\maketitle

\section*{Introduction}
For any topological group $G$, let $G^\delta$ denote $G$ with the same group structure, but considered as a discrete space.  The continuous homomorphism
$G^\delta\llra{\id}G$ induces a map of classifying spaces $BG^\delta\llra{i}BG$.  In \cite{miln83}, Milnor stated:

\noindent \textbf{Milnor's conjecture on the homology of Lie groups made discrete:} If $G$ is a Lie group, then $i$ induces isomorphisms\footnote{\ To ensure
readability, unadorned $H$ stands for the (co)homology of discrete groups throughout this paper, and `$\topo$' indicates singular (co)homology.}
\[
          H^n_\topo(BG,\Z/l)\ra H^n_\topo(BG^\delta,\Z/l)=H^n(G^\delta,\Z/l)
\]               \label{miln}
for any prime $l$.  The motivation for this (at first read, no doubt, surprising) conjecture was

\noindent \textbf{Friedlander's conjecture~\cite{frmi84}:} Let $k$ be an algebraically closed field and $l$ a prime distinct from the characteristic of $k$.
Let $G_k$ be an algebraic group over $k$. Then the natural map of group schemes $G(k)_k\ra G_k$ induces an isomorphism
\[
  H^n_\et(BG_k,\Z/l)\ra H^n_\et(BG_k(k),\Z/l)=H^n(G(k),\Z/l).
\]
Here $G(k)$ is the discrete group of $k$-rational points of the algebraic
group $G_k$, and the right-hand side is ordinary group cohomology.  On the
left, one has the \'etale cohomology of the simplicial scheme $BG_k$.
The construction of the comparison map between the two is akin to turning
a topological space into the discrete space made up of the set of its
points, save now one turns a scheme $G_k$ into $G(k)_k$, the coproduct of
copies of the terminal $k$-scheme $\spec(k)$, indexed by the
scheme-theoretic points of $G_k$.

The last chapter of Knudson~\cite{knud01} provides careful details and an excellent overview of the many partial results on these conjectures.  In addition,
the very end of section 1 of the present paper contains historical information on their origin that the reader is encouraged to become aware of.

Friedlander--Mislin~\cite{frmi84} proved Friedlander's conjecture over $k=\overline{\FF}_p$, the algebraic closure of a finite field.  The idea of this paper
is to exploit the interaction of the homologies of the discrete group $\tG(\overline{\FF}_p)$ and of $\tG(K)$ for large, algebraically closed fields $K$ to
attack more cases of Friedlander's conjecture.  Here $\tG$ is an integral form of a connected reductive algebraic group; such a $\tG$ is assumed to have been
fixed throughout.  Note that its group of complex points, $\tG(\CC)$, falls under the domain of both Milnor's and Friedlander's conjectures, and it is known
that for such groups the two are equivalent.

Our main result is formulated in terms of metric properties of the bar
complex for computing group homology.  Let $G$ be a discrete group and $R$
some ring of coefficients, on which $G$ is acting trivially.  Recall that
the bar complex is a functorial chain complex whose homology is
$H_*(G,R)$. The module of $n$-chains, $C_n(G)$, is the free $R$-module on
the basis set $G^n$; let $d_n$ denote the standard boundary map
$C_n(G)\xra{d_n}C_{n-1}(G)$, and $B_n(G)$ resp.\ $Z_n(G)$ the submodules
of $n$-boundaries and $n$-cycles.  Let the size $\|c\|$ of a chain $c\in
C_n$ mean the number of non-zero coefficients in the expression of $c$ as
formal linear combination of elements of $G^n$.  (If the coefficients $R$
were a normed abelian group, one would take the sum of the absolute values
of the coefficients, but throughout this paper we are concerned with prime
coefficients $R=\Z/l$.)

The \emph{filler norm} of a boundary $b\in B_n$ is defined as
\[
    \|b\|_\fil:=\min\big\{\|c\|\;\big|\;c\in C_{n+1}
                     \textup{\ such that\ }d_{n+1}(c)=b\big\}
\]
\begin{defn}   \label{ub}
$G$ satisfies a \emph{homological isoperimetric inequality for boundaries
in homological degree $n$ with coefficients $R$} if for all $K\in\N$,
\[
   \isop(K):=\sup\big\{\|b\|_\fil\;\big|\;b\in B_n
    \textup{\ such that\ }\|b\|=K\big\} < \infty
\]
\end{defn}
In words, the size of the shortest filler for a boundary $b$ can be
estimated from above in terms of the size of $b$ itself.  We call
$\isop$ the \emph{homological isoperimetric function} for $G$.

Fix now a $\tG$ as above, prime $p$, homological degree $n$ and
coefficients $\Z/l$, $l\neq p$.

{\bf Theorem A.} The following are equivalent:\\
$\bullet$ $\tG(\overline{\FF}_p)$ satisfies a homological isoperimetric inequality in degree $n$ with coefficients $\Z/l$\\
$\bullet$ Friedlander's conjecture holds for $H^n(\tG(k),\Z/l)$ for all algebraically closed fields $k$ of characteristic $p$.

The characteristic zero version involves, rather than an isoperimetric function of $\tG(\CC)^\delta$, an asymptotic isoperimetric function for
$\tG(\overline{\FF}_p)$ as $p$ ranges over the primes.

{\bf Theorem B.} Consider the statement:\\
(\asy) There exists a function $\asy:\N\ra\N$ with the property: for each $K\in\N$, for all sufficiently large primes $p$ (depending on $K$) one has that for
all $b\in B_n(\tG(\overline{\FF}_p))$ with $\|b\|=K$, $\|b\|_\fil\leq\asy(K)$.

\noindent $\bullet$ (\asy) implies that $H^n(\tG(\CC)^\delta,\Z/l)$ satisfies Milnor's conjecture; equivalently, that Friedlander's conjecture holds for
$H^n(\tG(k),\Z/l)$ for all algebraically closed fields $k$ of characteristic zero.\\
$\bullet$ The converse holds provided the homology of a maximal torus surjects on the homology of $\tG$; more precisely, if $\tG$ has a maximal torus $\tT$
(defined over the integers) such that for all  but finitely many primes $p$, the inclusion $\tT(\overline{\FF}_p)\inc\tG(\overline{\FF}_p)$ induces a
surjection
\[
      H_n(\tT(\overline{\FF}_p),\Z/l)\epi H_n(\tG(\overline{\FF}_p),\Z/l).
\]

This condition surfaces rather often in the study of Friedlander's conjecture, and is well-understood by a case-by-case analysis; perhaps it is enough to point
out that for any $\tG$, it holds for all primes $l$ that are large for $\tG$, i.e.\ that do not divide the order of the Weyl group of $\tG$ (and the list of
exceptional $l$ is typically much smaller).  See Section~\ref{isop} for more information on (removing) this obstacle from the converse implication.

Note that (\asy) neither implies that any particular
$\tG(\overline{\FF}_p)$ satisfies an isoperimetric inequality, nor is
implied by the existence of isoperimetric functions for individual
$\tG(\overline{\FF}_p)$ (unless those functions also happen to be
uniformly bounded in $p$).  (\asy) does imply that the uncountable group
$\tG(\CC)^\delta$ satisfies an isoperimetric inequality, but I do not know
the converse.

Let $k$ be an infinite field.  Set-theoretically, the cardinality of $H_n(\tG(k),\Z/l)$ is at most that of $k$.  Friedlander's conjecture predicts that for
algebraically closed $k$, $H_n(\tG(k),\Z/l)$ is isomorphic to the finite group $H_n^\topo(B\tG(\CC),\Z/l)$.  Since there is always a surjection, this is the
`smallest' value $H_n(\tG(k),\Z/l)$ can take.  The last result of this paper shows that as $k$ increases, the cardinality of $H_n(\tG(k),\Z/l)$ either grows as
fast as it can, or stays constant \emph{countable}.

{\bf Theorem C.} Fix $\tG$, $n$, $l$ and the characteristic $p\neq l$
through which our algebraically closed fields $k$ range ($p$ can be a
prime or zero).  One of the following two possibilities obtains:\\
$\bullet$ all the $H_n(\tG(k),\Z/l)$ are countable, all the groups
$\tG(k)$ possess the same isoperimetric function, and moreover every
extension $k\ra K$ between algebraically closed fields induces an
isomorphism $H_n(\tG(k),\Z/l)\llra{=}H_n(\tG(K),\Z/l)$,\\
$\bullet$ \emph{or} $H_n(\tG(k),\Z/l)$ has the cardinality of $k$ for all
uncountable $k$.

For $p>0$, thanks to Theorem A, the first alternative means the truth of
Friedlander's conjecture.  In particular, for an uncountable $k$ of
positive characteristic, $H_n(\tG(k),\Z/l)$ has either the value predicted
by Friedlander's conjecture, or the cardinality of $k$.  In characteristic
zero, Theorem C is far less useful; it does not restrict the values that
$H_n(\tG(\CC)^\delta,\Z/l)$ might take, and leaves open the possibility
that $\tG(\CC)^\delta$ possesses an isoperimetric function, though no
asymptotic isoperimetric function exists for the $\tG(\overline{\FF}_p)$.

\textsl{Terminological caveat.} Isoperimetric inequalities for boundaries
in the bar complex, at least with $\Z$ or $\RR$ coefficients, go back in
the literature to the 80's, prompted by Gromov's groundbreaking work on
bounded cohomology.  (See for example Matsumoto--Morita~\cite{morita1},
who refer to the condition `$\isop(K)\leq C\cdot K$ for some constant $C$'
as the \emph{uniform boundary property}.)  Isoperimetric functions for
$n$-balls in locally finite models of $K(G,1)$ appear under the name
\emph{higher Dehn functions}; see especially
Alonso--Wang--Pride~\cite{dehn}.  These generalize the classical
combinatorial Dehn function, or isoperimetric function, of finitely
presented groups.  Homological isoperimetric functions have also been
widely considered, especially in the context of hyperbolic groups, for
cycles (with $\Z$ or $\RR$ coefficients) on the universal cover of
suitable locally finite models of $K(G,1)$; see for example
Lang~\cite{lang00}.  It is not clear how the notions that pertain to
\emph{locally finite} models of $K(G,1)$ interact with isoperimetric
inequalities \emph{in the bar complex} --- not to mention that our groups,
such as $\tG(\overline{\FF}_p)$, are not finitely generated.  In this
paper, \emph{isoperimetric inequality} is always understood in the sense
of Def.~\ref{ub}.

\section{Heuristic}
The author discovered the relevance of isoperimetric inequalities by
analyzing a well-known bridge between the algebraic closures of finite
fields and uncountable algebraically closed fields.  Though the proofs can
be phrased without it, it is perhaps useful to give a blueprint of this
bridge, as the syntactic details of the argument may otherwise conceal the
simplicity of the main idea.

Let $\tP\subset\N$ be an infinite set of primes, and $\U$ any
non-principal ultrafilter on $\tP$.  It is an old observation that
\begin{equation}   \label{one}
     \prod_{\tP/\U} \overline{\FF}_p \iso \CC^\delta
\end{equation}
since both the ultraproduct on the left and the complex numbers (just as
an untopologized field) are algebraically closed fields of characteristic
zero, of the cardinality of the continuum.  The $\iso$ sign is to
emphasize how non-canonical the isomorphism is; it relies on Steinitz's
theorem, i.e.\ the possibility of a set-theoretic bijection between
transcendence bases (over the rationals) of the two sides.

Let $\tG$ be an algebraic group defined over the integers.  (\ref{one})
extends to give an isomorphism (non-canonically, and only as discrete
groups)
\begin{equation}   \label{two}
     \prod_{\tP/\U} \tG(\overline{\FF}_p) \iso \tG(\CC)^\delta
\end{equation}
The main ingredient in the proof of Friedlander's conjecture over $\overline{\FF}_p$ is the fact, proved earlier by Friedlander and (in special cases) by
Quillen, that for $p\neq l$
\[
   H_n(\tG(\overline{\FF}_p),\Z/l) \iso H_n^\topo(B\tG(\CC),\Z/l).
\]
Friedlander's conjecture asserts
\begin{equation}  \tag{?}
   H_n(\tG(\CC)^\delta,\Z/l) \iso H_n^\topo(B\tG(\CC),\Z/l).
\end{equation}
\emph{Supposing} that the functor $H_n(-,\Z/l)$ commutes with
ultraproducts,
\begin{equation}  \tag{?}
    \prod_{\tP/\U}H_n(\tG(\overline{\FF}_p),\Z/l) \iso
  H_n(\prod_{\tP/\U}\tG(\overline{\FF}_p),\Z/l)
\end{equation}
one would obtain Friedlander's conjecture:
\begin{equation}\begin{split} \label{three} \xymatrix{
 \overset{\phantom{blah}}{\prod_{\tP/\U}H_n(\tG(\overline{\FF}_p),\Z/l)}
       \ar@2{~}[d]\ar@2{~}[r]  &
 \overset{\phantom{blah}}{H_n(\prod_{\tP/\U}\tG(\overline{\FF}_p),\Z/l)}
       \ar@2{~}[d] \\
     H_n^\topo(B\tG(\CC),\Z/l) & H_n(\tG(\CC)^\delta,\Z/l)
}\end{split}\end{equation}
where the left-hand vertical isomorphism uses that, since the groups
$H_n^\topo(B\tG(\CC),\Z/l)$ are finite, canonically
\[
  \prod_{\tP/\U}H_n^\topo(B\tG(\CC),\Z/l)=H_n^\topo(B\tG(\CC),\Z/l).
\]
Ultraproducts, while slightly tamer, are nearly as badly behaved for
homological algebra as infinite products, and it is easy to see that the
functor $H_n(-,\Z/l)$ does not in general preserve them.  Nonetheless, one
has a natural comparison map
\[
     H_n(\prod_{\tP/\U}\tG(\overline{\FF}_p),\Z/l) \llra{[\iota]}
       \prod_{\tP/\U}H_n(\tG(\overline{\FF}_p),\Z/l)
\]
Most of the work goes into understanding the kernel and image of this homomorphism.  The reasons for falling back on the bar complex are its functoriality and
simple syntax, which make the interaction with ultraproducts much easier to analyze.  (That is also the reason for preferring to work with homology rather than
cohomology.)  The homomorphism $[\iota]$ turns out to be onto provided the homology of $\tG$ is supported on a maximal torus (and, I conjecture, in fact
always).  Via (\ref{three}), Friedlander's conjecture is seen to be equivalent to the injectivity of $[\iota]$.  The condition (\asy) results from a
combinatorial re-writing of this injectivity.

The positive characteristic case, Theorem A, is similar throughout but
much simpler; it uses ultrapowers of $\overline{\FF}_p$.  Theorem C
follows from the methods of the previous parts combined with an elementary
set-theoretic observation about constructible stratifications of algebraic
varieties over uncountable fields.

In homological degree $n=1$, the isoperimetric function of any group can
be understood completely in terms of its commutator width, and one can
establish the main properties of the isoperimetric functions of certain
groups (e.g.\ divisible abelian), in any homological degree, by hand.  It
seems to be challenging, however, to `reverse engineer' the deep and
beautiful work of Suslin~\cite{sus83}~\cite{sus84} in $K$-theory and
homological stability that yielded (in a range of dimensions) the
generalized isomorphism conjecture for $GL_N$ and $SL_N$, and to say
something about the isoperimetric functions associated to these groups;
not to mention, of course, establishing or refuting new cases of the
generalized isomorphism conjecture.  The difficulty is inherent, in part,
in the fact that our proof of Theorem B is non-constructive, i.e.\
proceeds by contradiction.  It is worth noting, however, the similarity
between Suslin's ``universal cycles'' and the isoperimetric condition in
the bar complex.  These threads will be pursued elsewhere; our goal here
is just to prove Theorems A, B and C.

\textsl{Historical remarks.} I am indebted to an anonymous referee for bringing the following to my attention.

The situation of Theorem B of this paper --- groups of the form $\tG(\CC)$ where $\tG$ is an integral form of a connected reductive algebraic group --- is
exactly what E.~Friedlander considered and discussed with many people during his stay in Princeton (1970--1975), when he first stated and investigated his
conjecture.  Milnor then generalized the conjecture to an arbitrary Lie group with finitely many connected components.  (Cf.\ the third sentence of
Milnor~\cite{miln83}: ``This paper is organized around the following conjecture which was suggested to the author by E.~Friedlander at least in the complex
case.'')  In the literature, the conjecture in this general form is called the ``Friedlander--Milnor Conjecture'' (or, as Milnor calls it, the ``Isomorphism
Conjecture'').

In other words, that part of Milnor's conjecture to which this paper has relevance is due to Friedlander.  The extension of the Isomorphism Conjecture to
algebraically closed base fields other than $\CC$, the ``Generalized Isomorphism Conjecture'', is due to Friedlander alone.

E.~Friedlander informs me that he and Charles Miller attempted to use ultraproducts to attack his conjecture in the 70's.

The first published use of ultraproducts is this context is due to Jardine~\cite{jard97}.  Other than the underlying idea of building uncountable algebraically
closed fields as ultraproducts of algebraic closures of finite fields, his methods and conclusions are distinct from ours.

\section{Ultraproducts of the bar complex}
The goal of this section is to construct a comparison homomorphism from
the homology of an ultraproduct of groups to the ultraproduct of their
homologies, and to give a necessary and sufficient condition for it to be
injective resp.\ surjective.  We only use elementary results on
ultraproducts and model theory in this paper, all contained in the
textbook Bell--Slomson~\cite{bellslo}.

The map is constructed via one particular device for computing group
homology, the bar complex.  For syntactic reasons, we spell out some
standard definitions in detail.  Fix a ring $R$ of coefficients, on which
all groups are understood to be acting trivially.  For a discrete group
$G$, the bar complex can be thought of as the simplicial homology of the
nerve of $G$ or, alternatively, as the result of tensoring with $-\tens_G
R$ the bar resolution of $R$ as trivial $G$-module.  The $n$-chains
$C_n(G)$ are the free $R$-module on the basis set $G^n$; we will write
basis elements as $\lc g_1,g_2,\dots,g_n\rc$ or $\vec{g}$.  The boundary
mapping $C_n\llra{d_n}C_{n-1}$ is defined on basis elements by
\begin{equation*}\begin{split}
\lc g_1,g_2,\dots,g_n\rc \mapsto \lc g_2,\dots,g_n\rc -
\lc g_1g_2,\dots,g_n\rc + \lc g_1,g_2g_3,\dots,g_n\rc - \dots    \\
+(-1)^{n-1}\lc g_1,g_2,\dots,g_{n-1}g_n\rc +
(-1)^n \lc g_1,g_2,\dots,g_{n-1}\rc
\end{split}\end{equation*}
(Here $d_1(\lc g_1\rc)=\lc\rc-\lc\rc=0\cdot\lc\rc$, where the empty tuple
$\lc\rc$ is the generator of $C_0(G)$, and $C_{-1}(G)$ is by definition
zero.)  If $z$ is a cycle, we write $[z]$ for the homology class it
represents.

Let $G_\lambda$, $\lambda\in\Lambda$, be a set of discrete groups.  Let
$\U$ be an ultrafilter on $\Lambda$, and let $\G$ denote the corresponding
ultraproduct $\prod_{\Lambda/\U} G_\lambda$. If $\phi(\lambda)$ is a
mathematical statement containing the parameter $\lambda$ ranging over
$\Lambda$, we will abbreviate as
\[   \U \models\;\phi    \]
the statement ``the set of $\lambda\in\Lambda$ for which $\phi(\lambda)$
is true, belongs to $\U$''.  (Remark: only the variable $\lambda$ will be
used in this role.  Though the notation is suggestive, it is meant to be
just a typographical device.  In particular, $\phi$ will be typically
phrased in the meta-language, and is not necessarily assumed to be
equivalent to a first-order formula in the language of rings and groups.)

Consider now the ultraproduct (as $R$-modules)  $\prod_{\Lambda/\U}
C_n(G_\lambda)$.  One has an $R$-linear map
\[
   C_n(\G)\llra{\iota}\prod_{\Lambda/\U} C_n(G_\lambda)
\]
defined on basis elements as follows: given $\vec{g}\in\G^n$, choose
representatives $\{\vec{g}_\lambda\in G_\lambda^n\,|\,\lambda\in\Lambda\}$
for it; the collection
\[
   \big\{1\cdot\vec{g}_\lambda\,|\,\lambda\in\Lambda\big\}
\]
gives a well-defined element $\iota(\vec{g})$ of the ultraproduct
$\prod_{\Lambda/\U} C_n(G_\lambda)$.

\begin{prop}
$\iota$ is injective.
\end{prop}

\begin{proof}
Let $X$ be an arbitrary set, and fix a $k$-tuple
$\lc r_1,r_2,\dots,r_k\rc$ of elements of $R$.  To say that the formal
expression $r_1\cdot x_1+r_2\cdot x_2+\dots+r_k\cdot x_k$ (where the $x_i$
are thought of as variable, ranging over $X$) equals $0$ in the free
$R$-module with basis $X$ amounts to a first-order formula
\[
\bigvee_{\substack{I_1\sqcup I_2\sqcup\dots\sqcup I_p\\=\{1,2,\dots,k\}}}
\;\bigwedge_{q=1}^p\;\bigwedge_{i,j\in I_q} x_i=x_j
\]
where the disjunction is over all partitions of $\{1,2,\dots,k\}$ into
subsets $\{I_1,I_2,\dots,I_p\}$ in such a way that $\sum_{i\in I_q} r_i=0$
for each $q=1,2,\dots,p$.  Call this formula $\theta(x_1,x_2,\dots,x_k)$.
Given $c=\sum r_i\,\vec{g}_i$ in $C_n(\G)$ and representatives
$\{\vec{g}_{i,\lambda}\in G_\lambda^n\,|\,\lambda\in\Lambda\}$ for
$\vec{g}_i$, the following are equivalent:
\begin{align*}
& \iota(c)=0\textup{ in\ }\prod_{\Lambda/\U} C_n(G_\lambda) \\
& \dlra \U\models\quad \sum r_i\,\vec{g}_{i,\lambda}=0
                   \textup{ in\ }C_n(G_\lambda)             \\
& \dlra \U\models\quad
\theta(\vec{g}_{1,\lambda},\vec{g}_{2,\lambda},\dots,\vec{g}_{k,\lambda})\\
& \dlra \theta(\vec{g}_1,\vec{g}_2,\dots,\vec{g}_k)
                   \textup{ holds in\ }C_n(\G)              \\
& \dlra c=\sum r_i\,\vec{g}_i=0\textup{ in\ }C_n(\G)
\end{align*}
by \Los's theorem.
\end{proof}

$\prod_{\Lambda/\U} C_n(G_\lambda)$ can be equipped with a `boundary' map
$\hat{d}_n: \prod_{\Lambda/\U} C_n(G_\lambda)\ra
\prod_{\Lambda/\U} C_{n-1}(G_\lambda)$, which is simply the ultraproduct
of the boundary mappings connecting the individual $C_n(G_\lambda)$.  It
is therefore linear and one checks (via basis elements) that
$\iota(d_n(c))=\hat{d}_n(\iota(c))$ for
$c\in\prod_{\Lambda/\U} C_n(G_\lambda)$.

Now for each homology class in $H_n(\G,R)$ take a representing cycle
$z=\sum r_i\,\vec{g}_i$ from $Z_n(\G)$.  Choosing representatives
$\{\vec{g}_{i,\lambda}\in G_\lambda^n\,|\,\lambda\in\Lambda\}$ for each
$\vec{g}_i$, one has that
\[
  \U\models\quad \sum r_i\,\vec{g}_{i,\lambda}\,\in Z_n(G_\lambda)
\]
hence $\iota(z)$ represents an element in (the ultraproduct as
$R$-modules) $\prod_{\Lambda/\U} H_n(G_\lambda,R)$.  This
element is independent of the choice of representative $z$ taken
in its homology class.  Indeed, if $z'$ is another such, then
$z-z'=d_{n+1}(c)$ for some $c\in C_{n+1}(\G)$, implying
\[
  \U\models\quad \iota(z)\textup{ and\ }\iota(z')\textup{ are homologous
   in\ }Z_n(G_\lambda).
\]
One therefore obtains a map
\begin{equation}   \label{comp}
    H_n(\G,R) \llra{[\iota]} \prod_{\Lambda/\U} H_n(G_\lambda,R)
\end{equation}
\begin{remark}
It is also true that $\hat{d}_n\circ\hat{d}_{n+1}=0$ and
$\ker\hat{d}_n/\im\hat{d}_{n+1}$ is canonically isomorphic to
$\prod_{\Lambda/\U} H_n(G_\lambda,R)$, but we won't need this.
\end{remark}

Using the fact that the algebraic structure on homology classes is definable directly on cycle representatives, one sees that $[\iota]$ is $R$-linear.  We wish
to understand the kernel and image of $[\iota]$.  This turns out to be more tedious for the case of an infinite $R$, nor does that case have relevance to
Friedlander's conjecture.  (See the Appendix of Milnor~\cite{miln83} for an investigation of $H_*^\topo(BG^\delta)$ with rational or real coefficients.)
Henceforth we \emph{assume the cardinality of $R$ to be finite}, and introduce a partial inverse to $\iota$.

Recall that the \emph{size} $\|c\|$ of an element of a free $R$-module with specified basis is the number of basis elements occurring with non-zero
coefficients.

\begin{defn}
An element $\hat{c}$ of $\prod_{\Lambda/\U} C_n(G_\lambda)$ is said to be
\emph{$\U$-uniformly bounded} (or simply \emph{bounded}) if there exists
$K<\infty$ such that for some (equivalently, all) representatives
$\{c_\lambda\in C_n(G_\lambda)\,|\,\lambda\in\Lambda\}$ of $\hat{c}$,
\[
   \U\models\; \|c_\lambda\|\leq K
\]
\end{defn}

It is immediate that for any $c\in C_n(\G)$, $\iota(c)$ is bounded; bounded chains form a submodule of $\prod_{\Lambda/\U} C_n(G_\lambda)$; and
$\hat{d}(\hat{c})$ is bounded if $\hat{c}$ is so.

Let $\hat{c}\in\prod_{\Lambda/\U} C_n(G_\lambda)$ be bounded in size by
$K$.  Choose representatives $\{c_\lambda\in
C_n(G_\lambda)\,|\,\lambda\in\Lambda\}$ of $\hat{c}$, and write each
$c_\lambda$ with $|c_\lambda|\leq K$ in \emph{some} way as an
\emph{ordered} sum of basis elements,
$c_\lambda=\sum_{i=1}^{l_\lambda} r_i\,\vec{g}_{i,\lambda}$,
$l_\lambda\leq K$.  This allows one to define a map recording
`coordinates'
\[
  \{c_\lambda\in C_n(G_\lambda)\textup{\ such that\ }|c_\lambda|\leq K\}
          \xra{\coo}R^{\leq K}   \]
(where $R^{\leq K}$ is the set of ordered tuples from $R$ of
size at most $K$) by sending $c_\lambda=\sum_{i=1}^{l_\lambda}
r_i\,\vec{g}_{i,\lambda}$ to $\coo(c_\lambda)\overset{\de}{=}
\lc r_1,r_2,\dots,r_{l_\lambda}\rc$. The map $f$ from the appropriate
element of $\U$ to $R^{\leq K}$ defined by
$\lambda\mapsto\coo(c_\lambda)$ partitions a member of the
ultrafilter into finitely many disjoint subsets via $f^{-1}(t)$ as
$t$ ranges over the elements of $R^{\leq K}$.  So for exactly
one tuple $t_0$ will the set $f^{-1}(t_0)$ belong to $\U$.  Let
that $t_0$ be $\lc r_1,r_2,\dots,r_l\rc$ and write $U$ for
$f^{-1}(t_0)$; then one has that for all $\lambda\in U$,
\[
   c_\lambda=\sum_{i=1}^l r_i\,\vec{g}_{i,\lambda}
\]
for well-defined $r_i\in R$, $\vec{g}_{i,\lambda}\in
G_\lambda^n$. For a fixed $i$, the collection
$\{\vec{g}_{i,\lambda}\;|\;\lambda\in U\}$ (extended by arbitrary
$\vec{g}_{i,\lambda}$ for $\lambda\in\Lambda\setminus U$, if
necessary) can be thought of as an element $\vec{g}_i$ in the
ultraproduct $\G^n$.  Introduce the notation
\[   \tau(\hat{c})\,\overset{\de}{=}\,
       \sum_{i=1}^l a_i\,\vec{g}_i\,\in\,C_n(\G)  .  \]

Obviously $\iota(\tau(\hat{c}))=\hat{c}$.  But since $\iota$ is
injective, $\iota$ and $\tau$ must be inverse bijections between
$C_n(\G)$ and the submodule of $\prod_{\Lambda/\U} C_n(G_\lambda)$
consisting of bounded chains.  In particular, $\tau$ is
independent of the choices made, $R$-linear, and
$d_n(\tau(\hat{c}))=\tau(\hat{d}_n(\hat{c}))$ for any bounded
$\hat{c}\in\prod_{\Lambda/\U} C_n(G_\lambda)$.

\begin{cor}
Call a class in $\prod_{\Lambda/\U} H_n(G_\lambda,R)$
\emph{bounded} if it can be represented by a bounded cycle in
$\prod_{\Lambda/\U}C_n(G_\lambda)$.  Bounded classes form a
submodule of $\prod_{\Lambda/\U} H_n(G_\lambda,R)$, which equals
the image of $[\iota]$.
\end{cor}

\begin{cor}  \label{image}
$[\iota]$ is surjective if and only if there exists $K<\infty$ such that
for $\U$-most $\lambda$, every homology class in $H_n(G_\lambda,R)$
contains a cycle of size at most $K$.
\end{cor}

\begin{prop}  \label{iff}
$[\iota]$ is injective if and only if the following holds:

$(\bigstar)$ For any bounded $\hat{b}\in\prod_{\Lambda/\U}
C_n(G_\lambda)$, if $\U\models b_\lambda\in B_n(G_\lambda)$, then
there exists a bounded $\hat{c}\in\prod_{\Lambda/\U}
C_{n+1}(G_\lambda)$ such that $\U\models\;
b_\lambda=d_{n+1}(c_\lambda)$.
\end{prop}

\begin{proof}
Suppose $(\bigstar)$ holds.  Let $b\in Z_n(\G)$ be a cycle
representative of a homology class in $H_n(\G,R)$ that is sent
to zero by $[\iota]$.  That is to say, $\U\models \iota(b)\in
B_n(G_\lambda)$.  $\iota(b)$ is bounded, hence by assumption there
exists a bounded $\hat{c}\in\prod_{\Lambda/\U} C_{n+1}(G_\lambda)$
such that $\hat{d}_{n+1}(\hat{c})=\iota(b)$.  But then
$d_{n+1}(\tau(\hat{c}))=\tau(\hat{d}_{n+1}(\hat{c}))=\tau(\iota(b))=b$.
Thence $b$ is a boundary, so represents the zero homology class.

Conversely, assume $[\iota]$ is injective, and let $\hat{b}$ be bounded
and $\U$-almost everywhere a boundary.  Then
$d_n(\tau(\hat{b}))=\tau(\hat{d}_n(\hat{b}))=\tau(0)=0$, so
$\tau(\hat{b})$ represents a homology class in $H_n(\G,R)$.
$[\iota]([\tau(\hat{b})])=0$ since $\iota(\tau(\hat{b}))=\hat{b}$ is
$\U$-almost everywhere a boundary by assumption.  By the injectivity of
$[\iota]$, there must exist $c\in C_{n+1}(\G)$ such that
$d_{n+1}(c)=\tau(\hat{b})$.  $\iota(c)$ is bounded and
\[
   \hat{d}_{n+1}(\iota(c))=\iota(d_{n+1}(c))=\iota(\tau(\hat{b}))=\hat{b}
\]
so $(\bigstar)$ is satisfied.
\end{proof}

In effect, this says that $(\bigstar)$ holds if and only if
$[\hat{z}]\mapsto[\tau(\hat{z})]$ is the inverse bijection to
$[z]\mapsto[\iota(z)]$ between $H_n(\G,R)$ and the module of bounded
homology classes.

\section{From Friedlander's conjecture to isoperimetric functions}   \label{isop}
The key result that makes the previous section applicable to Friedlander's conjecture is due to Friedlander~\cite{friedlander}~\cite{frmi84}, based on work of
Quillen~\cite{quil70a}:

\begin{thm} Let $p\neq l$ be primes.   \label{quil}
$H_n(\tG(\overline{\FF}_p),\Z/l)$ is isomorphic to
$H_n^\topo(B\tG(\CC),\Z/l)$.
\end{thm}

This latter group is finite (and known, as a function of $\tG$, $n$ and
$l$).  Let us introduce the notation $|H_{\tG,n,l}|$ for the common value
of the cardinalities of the groups $H_n(\tG(\overline{\FF}_p),\Z/l)$,
$l\neq p$, and $H_n^\topo(B\tG(\CC),\Z/l)$.

\begin{lemma}    \label{1st}
$\card H_n(\tG(\CC)^\delta,\Z/l)\geq|H_{\tG,n,l}|$.  Friedlander's conjecture holds for the Lie group $\tG(\CC)$ if and only if $\card
H_n(\tG(\CC)^\delta,\Z/l)=|H_{\tG,n,l}|$.
\end{lemma}

\begin{proof}
Both statements follow from the theorem of Milnor~\cite{miln83} that
$H_n(G^\delta,\Z/l)\llra{i}H_n^\topo(BG,\Z/l)$ is surjective for any Lie
group $G$ with finitely many components.
\end{proof}

The next observation (also well-known) contains Lemma~\ref{1st}; they are
repeated just to emphasize the parallel.  Note that Milnor's proof is
purely topological and applies to real Lie groups as well.

\begin{lemma}    \label{2nd}
Let $k$ be an algebraically closed field of characteristic $p\neq l$.
Then $\card H_n(\tG(k),\Z/l)\geq|H_{\tG,n,l}|$.  Friedlander's conjecture
holds for $\tG_k$ if and only if $\card H_n(\tG(k),\Z/l)=|H_{\tG,n,l}|$.
\end{lemma}

\begin{proof}
By a theorem of Friedlander and Mislin, the map
$H^n_\et(B\tG_k,\Z/l) \ra H^n(\tG(k),\Z/l)$ concerned in the generalized
isomorphism conjecture is injective.  Now
$H^n_\et(B\tG_k,\Z/l)\iso H^n_\et(B\tG_{\overline{\FF}_p},\Z/l)\iso
H^n(\tG(\overline{\FF}_p),\Z/l)\iso H_n(\tG(\overline{\FF}_p),\Z/l)$ by
virtue of the invariance of \'etale cohomology under algebraically closed
field extensions and the truth of the generalized isomorphism conjecture
over $\overline{\FF}_p$; and if $H^n(\tG(k),\Z/l)$ or $H_n(\tG(k),\Z/l)$
is finite, then they are isomorphic.  These facts imply both parts.
\end{proof}

We are now ready to prove one direction of Theorem A:

\begin{prop}
If Friedlander's generalized isomorphism conjecture holds for
$H_n(\tG(k),\Z/l)$ for all algebraically closed fields $k$ of
characteristic $p$, then $\tG(\overline{\FF}_p)$ satisfies an
isoperimetric inequality in homological degree $n$ with coefficients
$\Z/l$.
\end{prop}

\begin{proof}
Consider any non-principal ultrafilter $\U$ on any countable set
$\Lambda$, and write $P$ for the ultrapower
$\prod_{\Lambda/\U}\overline{\FF}_p$.  $P$ is an algebraically closed
field of characteristic $p$, of the cardinality of the continuum.  As
$\tG$ is first-order definable in the language of rings,
$\prod_{\Lambda/\U}\tG(\overline{\FF}_p)$ is canonically isomorphic to
$\tG(P)$.  Apply the comparison homomorphism (\ref{comp}):
\begin{equation}   \label{inj2}
    H_n(\tG(P),\Z/l)\llra{[\iota]}
    \prod_{\Lambda/\U} H_n(\tG(\overline{\FF}_p),\Z/l)
\end{equation}
Since $H_n(\tG(\overline{\FF}_p),\Z/l)$ is finite, the right-hand side is
isomorphic to $H_n(\tG(\overline{\FF}_p),\Z/l)$.  $[\iota]$ is surjective,
with a splitting induced by the inclusion $\overline{\FF}_p\inc P$ (which
can also be thought of as the diagonal
$\overline{\FF}_p\ra\prod_{\Lambda/\U}\overline{\FF}_p$); see
Cor.~\ref{image}.

Proceed by contradiction.  Fix some size $K\in\N$, and let $B_K$ denote
the set of boundaries $b\in B_n(\tG(\overline{\FF}_p))$ with $\|b\|=K$.
Suppose
\[
    \big\{\|b\|_\fil\;\big|\;b\in B_K\big\} \subseteq\N
\]
were unbounded.  Then there would exist an infinite subset
$\Lambda\subseteq B_K$ such that for any infinite $U\subseteq\Lambda$,
\[
    \big\{\|b\|_\fil\;\big|\;b\in U\big\}
\]
is still unbounded.  This $\Lambda$ will serve as the index set for a
comparison map of the type (\ref{inj2}); the non-principal ultrafilter
$\U$ on $\Lambda$ can be arbitrary.

One has a tautologous element $\hat{s}\in\prod_{\Lambda/\U}
B_n(G_\lambda)$, associating to $\lambda\in\Lambda$ itself, since
$\Lambda\subseteq B_n(\tG(\overline{\FF}_p))$.  By construction, $\hat{s}$
is $\U$-bounded in size by $K$.  The existence of a $\U$-bounded
$\hat{c}\in\prod_{\Lambda/\U} C_{n+1}(G_\lambda)$ such that
$\U\models\;s_\lambda=d_{n+1}(c_\lambda)$ would mean that there exist a
constant $K_1<\infty$ and $U\in\U$ such that for all $b\in U$,
$\|b\|_\fil\leq K_1$.  Since $\U$ is non-principal, $U$ would be an infinite
subset of $\Lambda$, contradicting the choice of $\Lambda$.  Therefore
property ($\bigstar$) of Cor.~\ref{iff} fails, and $[\iota]$ cannot be
injective.

But that means that the cardinality of $H_n(\tG(P),\Z/l)$ exceeds that of
$H_n(\tG(\overline{\FF}_p),\Z/l)$, and so (cf.\ Lemma~\ref{2nd}) the
generalized isomorphism conjecture fails over $P$, the (unique)
algebraically closed field of characteristic $p$ that has the cardinality
of the continuum.
\end{proof}

The other direction of Theorem A follows easily from the comparison map
(\ref{inj2}), injectivity condition ($\bigstar$), and a theorem of
Friedlander-Mislin stating that if the generalized isomorphism conjecture
holds for one algebraically closed field of infinite transcendence degree
over its prime field, then it holds for all algebraically closed fields
within that characteristic.  However, we prefer to give a completely
elementary and self-contained proof of that direction in
Section~\ref{algclo}.

Let us turn to the hard part of Theorem B.  Let $\tP$ be any infinite set
of primes, and $\U$ any non-principal ultrafilter on $\tP$.  The
ultraproduct of discrete groups $\prod_{\tP/\U}\tG(\overline{\FF}_p)$ is
(non-canonically) isomorphic to $\tG(\CC^\delta)=\tG(\CC)^\delta$, a
complex algebraic Lie group made discrete.  Apply the comparison
homomorphism (\ref{comp}) to the family $\tG(\overline{\FF}_p)$,
$p\in\tP$:
\begin{equation}   \label{inj}
    H_n(\prod_{\tP/\U}\tG(\overline{\FF}_p),\Z/l)\llra{[\iota]}
    \prod_{\tP/\U} H_n(\tG(\overline{\FF}_p),\Z/l)
\end{equation}

For all $p\neq l$, a fortiori for all but finitely many $p\in\tP$, $H_n(\tG(\overline{\FF}_p),\Z/l)$ is isomorphic to the finite group
$H_n^\topo(B\tG(\CC),\Z/l)$, so the right-hand side of (\ref{inj}) is isomorphic to $H_n^\topo(B\tG(\CC),\Z/l)$.  Though the argument is analogous to
characteristic $p$, a point has to be overcome in order to deduce the condition (\asy) from Friedlander's conjecture.

\begin{lemma}   \label{onto}
Suppose $\tG$ has a maximal torus $\tT$ (defined over the integers) such
that for all but finitely many primes $p$, the inclusion
$\tT(\overline{\FF}_p)\inc\tG(\overline{\FF}_p)$ induces a surjection
\[
  H_n(\tT(\overline{\FF}_p),\Z/l)\epi H_n(\tG(\overline{\FF}_p),\Z/l).
\]
Then the $[\iota]$ of \textup{(\ref{inj})} is surjective, for any
non-principal ultrafilter $\U$ on any infinite set of primes $\tP$.
\end{lemma}

\begin{proof}
We wish to prove the following: there exists a bound $f_\tG(n,l)<\infty$
such that for all but finitely many primes $p$, each homology class
$\alpha\in H_n(\tG(\overline{\FF}_p),\Z/l)$ has a cycle representative
$z_{p,\alpha}\in Z_n(\tG(\overline{\FF}_p))$ with
$\|z_{p,\alpha}\|\leq f_\tG(n,l)$.  By Cor.~\ref{image}, this implies (and
is in fact equivalent to) the conclusion.

Such a bound exists for the torus $\tT=GL_1\times\dots\times GL_1$ of rank
$r$.  Assume $p\neq l$.  Identifying the $l^\th$-power roots of unity in
$\overline{\FF}_p$ with $\Z/l^\infty$, one gets an injection
$\Z/l^\infty\ra\overline{\FF}_p^\times$ such that
$(\Z/l^\infty)^r\ra\tT(\overline{\FF}_p)$ induces isomorphism on
$H_*(-,\Z/l)$.  One can take
$f_\tT(n,l)=\max\{\|z_1\|,\|z_2\|,\dots,\|z_N\|\}$ where the cycles $z_i$
span the (finite) group $H_n((\Z/l^\infty)^r,\Z/l)$.

Under the assumption that the homology of a maximal torus surjects on the
homology of $\tG$, one can take $f_\tG(n,l)=f_\tT(n,l)$ for the torus
$\tT$ of the same rank as $\tG$.
\end{proof}

\noindent\textbf{Discussion.}
For the sake of completeness, let us recall how `cheap' this assumption
is.  There are several well-known ways to investigate
$H_n(\tT(\overline{\FF}_p),\Z/l)\epi H_n(\tG(\overline{\FF}_p),\Z/l)$.

The functoriality of Friedlander's isomorphism
$H_n(\tG(\overline{\FF}_p),\Z/l)\iso H_n^\topo(B\tG(\CC),\Z/l)$ is subtle,
as it depends on an embedding of the Witt vectors of $\overline{\FF}_p$ in
$\CC$.  However, by making choices simultaneously for $\tG$ and its split
maximal torus $\tT$, one obtains a commutative diagram
\[\xymatrix{
    H_n(\tT(\overline{\FF}_p),\Z/l)\ar[r]\ar[d]^{\iso} &
                     H_n(\tG(\overline{\FF}_p),\Z/l)\ar[d]^{\iso} \\
    H_n^\topo(B\tT(\CC),\Z/l)\ar[r] &  H_n^\topo(B\tG(\CC),\Z/l)
}\]
On the topological side, one has a surjection
$H_n^\topo(B\tT(\CC),\Z/l)\epi H_n^\topo(B\tG(\CC),\Z/l)$ when (for
example) $l$ is prime to the order of the Weyl group of $\tG$; one way to
see this is to approximate the classifying space of a Lie group by
manifolds, and use Becker--Gottlieb transfer.  See
Feshbach~\cite{feshbach}.

For Chevalley groups $\tG$, one can also argue purely group-theoretically.
Suppose $l\neq p$, $l\nmid|\tW|$, and let the $p$-power $q$ be such that
$\FF_q$ contains $l^\th$ roots of unity.  By a theorem of
Chevalley~\cite{cheval}, there exists a split maximal torus $\tT$ of $\tG$
such that $\tT(\FF_q)$ contains a Sylow $l$-subgroup of $\tG(\FF_q)$.
(This is because $\big[\tG(\FF_q):\tT(\FF_q)\big]$ will be prime to $l$.)
Therefore $\tT(\FF_q)\inc\tG(\FF_q)$ induces a surjection
\[
   H_n(\tT(\FF_q),\Z/l)=H_n(\textup{Syl}_l\tT(\FF_q),\Z/l)
    =H_n(\textup{Syl}_l\tG(\FF_q),\Z/l)\epi H_n(\tG(\FF_q),\Z/l).
\]
Let $\FF_q$ be cofinal in $\overline{\FF}_p$ such that
$q\equiv 1\pmod{l}$.  Since the tori $\tT$ can be chosen compatibly,
there results a surjection
$H_n(\tT(\overline{\FF}_p),\Z/l)\epi H_n(\tG(\overline{\FF}_p),\Z/l)$.

If $l$ is a torsion prime for $\tG$, one need not have a surjection
$H_n(\tT(\overline{\FF}_p),\Z/l)\epi H_n(\tG(\overline{\FF}_p),\Z/l)$.
Nonetheless, a Sylow subgroup of $\tG(\FF_q)$ is always contained in the
normalizer of a torus, which is an extension of a torus by the Weyl group.
By making compatible choices as $\FF_q$ increases to $\overline{\FF}_p$,
one obtains a surjection
$H_n(N_\tT(\overline{\FF}_p),\Z/l)\epi H_n(\tG(\overline{\FF}_p),\Z/l)$
with a short exact sequence
$\te\ra\tT(\overline{\FF}_p)\ra N_\tT(\overline{\FF}_p)\ra\tW\ra\te$.

The Lyndon-Hochschild-Serre spectral sequence for this extension has the
form
\[
   H_i\big(\tW,H_j(\tT(\overline{\FF}_p),\Z/l)\big)\dfra
                           H_{i+j}(N_\tT(\overline{\FF}_p),\Z/l).
\]
Since all homology groups involved are finite, the spectral sequence converges (in any total degree) in a finite number of steps.  In principle at least, one
can check that the homology of $\tW$ has cycle representatives (with twisted coefficients) whose size is bounded independently of $p$; analyzing the
differentials in the spectral sequence, presumably so does $H_n(N_\tT(\overline{\FF}_p),\Z/l)$ and, eventually, $H_n(\tG(\overline{\FF}_p),\Z/l)$ --- for all
$\tG$, $n$, $p$ and $l$.  At this stage, it does not seem worthwhile to spell out these details for exceptional $l$.  (Recall that the implication from the
existence of an asymptotic isoperimetric function to the truth of Friedlander's conjecture, to be proved in section~\ref{algclo}, holds unconditionally.)

We return to the proof of Theorem B now.

\begin{prop}   \label{hardb}
Let $\tG$, $n$, $l$ be such that the conclusion of Lemma~\ref{onto} holds.  Friedlander's conjecture for $H_n(\tG(\CC)^\delta,\Z/l)$ implies that the
asymptotic isoperimetric function of Theorem B exists.
\end{prop}

\begin{proof}
Fix some $K\in\N$.  Write $B_K(p)$ for the set of boundaries
$b\in B_n(\tG(\overline{\FF}_p))$ with $\|b\|=K$.  By contradiction,
assume: for all $K_1\in\N$, there exist infinitely many primes $p$ such
that $\sup\big\{\|b\|_\fil\;|\;b\in B_K(p)\big\}>K_1$.

That would allow one to find an infinite set $\tP=\{p_0,p_1,p_2,\dots\}$ of primes and for each $p\in\tP$ some boundary $b_p\in B_K(p)$ with the property that
for any infinite subset $U\subseteq\tP$, the set $\{\|b_p\|_\fil\;|\;p\in U\}$ is unbounded.  (Let $p_0$ and $b_{p_0}$ be arbitrary, and having found $p_n$,
pick $p_{n+1}\not\in\{p_0,p_1,\dots,p_n\}$ and $b_{p_{n+1}}\in B_K(p_{n+1})$ such that $\|b_{p_{n+1}}\|_\fil>\|b_{p_n}\|_\fil$.)  Let this $\tP$ be the
infinite set of primes with which the comparison homomorphism (\ref{inj}) is constructed.  If Friedlander's conjecture holds, then by Lemma~\ref{1st}, the two
sides of (\ref{inj}) have the same finite cardinality.  So if $[\iota]$ is surjective, it must be injective too, and condition $(\bigstar)$ of Cor.~\ref{iff}
must be satisfied.  On the other hand, for the $\U$-bounded element $\{p\mapsto b_p\}\in\prod_{\tP/\U} B_n(\tG(\overline{\FF}_p))$ there cannot exist a
$\U$-bounded $\hat{c}\in\prod_{\tP/\U} C_{n+1}(\tG(\overline{\FF}_p))$ such that $\U\models\;b_p=d_{n+1}(c_p)$: since $\U$ is non-principal, that would mean
that for some infinite subset $U\subseteq\tP$, there does exist $K_1$ such that for all $p\in U$, $\|b_p\|_\fil<K_1$, contradicting the choice of $\tP$.

Therefore, under Friedlander's conjecture, one can find $\asy(K)=K_1<\infty$ such that for all but finitely many primes $p$,\, $\sup\big\{\|b\|_\fil\;|\;b\in
B_K(p)\big\}\leq K_1$.
\end{proof}

\section{From isoperimetric functions to Friedlander's conjecture}    \label{algclo}
One can phrase the mathematics behind the other directions of Theorems A
and B in two ways, different only linguistically.  One is the language of
constructible subsets of varieties over algebraically closed fields,
Chevalley's theorem on the image of constructible sets under regular maps
being constructible, base extensions between algebraically closed fields,
and specialization (this is the spirit of the next section) and the other
is the language of sets definable in the first-order theory of
algebraically closed fields, Tarski's theorem on quantifier elimination,
and the first-order Lefschetz principle.  Considering the syntax of the
statements involved, the second approach seems much more convenient, and
that is what we will use.

\textbf{Conventions.}  The algebraic group $\tG$ defined over the
integers, homological degree $n$, and finite ring of coefficients $R=Z/l$
will be fixed once and for all.  Variables will range over the
algebraically closed field $k$; that makes the group of $k$-rational
points $\tG(k)$ and the group operations on $\tG(k)$ first-order
expressible in the language of rings.  Observe that none of ``chain'',
``cycle'' and ``boundary'' in the bar complex are first-order expressible.
However, \emph{for any given choice of the bounds $K$, $K_1$}, each of
\begin{center}
`` chain $c\in C_n(\tG(k))$ with $\|c\|=K$ ''\\
`` cycle $z\in Z_n(\tG(k))$ with $\|z\|=K$ ''\\
`` boundary $b\in B_n(\tG(k))$ with $\|b\|=K$ and $\|b\|_\fil=K_1$ ''
\end{center}
is first-order expressible.  (Code a chain $c\in C_n$ of size $K$ as
$K\cdot l$ many $n$-tuples of elements of $\tG(k)$, exploit the
first-order definition of the bar differential $d_n$ and the fact that the
equality of two expressions that are unordered formal $R$-linear
combinations is first-order.)

For any $K,K_1,K_2\in\N$, consider the sentence $\Phi_{K,K_1,K_2}$

\noindent
`` For every $b\in C_n(\tG(k))$ with $\|b\|=K$, if there exists
$u\in C_{n+1}(\tG(k))$ with $\|u\|=K_1$ such that $d_{n+1}(u)=b$, then
there exists $c\in C_{n+1}(\tG(k))$ with $|c|\leq K_2$ such that
$d_{n+1}(c)=b$. ''

By Tarski's theorem, $\Phi_{K,K_1,K_2}$ either holds in all algebraically
closed fields $k$ of a given characteristic, or none.  But the countable
conjunction $\bigwedge_{K_1\in\N}\Phi_{K,K_1,K_2}$ means precisely
\[
    \textup{For any boundary with\ }\|b\|=K,
              \textup{\ one has\ }\|b\|_\fil\leq K_2.
\]
\begin{cor}  \label{same}
Fix any $K\in\N$.  As $k$ ranges through algebraically closed fields
within any given characteristic,
\[ \sup\big\{
      \|b\|_\fil\;|\;b\in B_n(\tG(k))\textup{\ such that\ }\|b\|=K
       \big\}  \in\N\cup\infty
\]
stays the same.
\end{cor}

So it makes sense to talk of ``the isoperimetric function of $\tG$ in
characteristic $p$'', $p$ a prime or zero, provided this supremum is
finite for all $K\in\N$.

\begin{prop}   \label{oneway}
If $\tG$ satisfies an isoperimetric inequality in characteristic $p$ ($p$
a prime or zero) and $H_n(\tG(k),R)$ is finite for one particular
algebraically closed $k$ of that characteristic, then the groups
$H_n(\tG(k),R)$ are isomorphic for all algebraically closed $k$ of
characteristic $p$.
\end{prop}

\begin{proof}
Let $k_1$ be such that the cardinality of $H_n(\tG(k),R)$, as $k$ varies
over algebraically closed fields in characteristic $p$, is minimal at
$k=k_1$.  Write $|H|$ for that least cardinality; by assumption,
$|H|<\infty$.  For any $K\in\N$, set
$K_1=\max\{\isop(1),\isop(2),\dots,\isop(2K)\}$ and consider the sentence

\noindent
$\big(\Psi_K\big)$\quad `` Given $z_i\in Z_n(\tG(k))$ with
$\|z_i\|\leq K$, $i=1,2,\dots,|H|+1$, there exist
$1\leq i\neq j\leq |H|+1$ and $c\in C_{n+1}(\tG(k))$ with $|c|\leq K_1$
such that $z_i-z_j=d_{n+1}(c)$. ''

\noindent
Since this is first-order and holds over $k=k_1$, it holds for all
algebraically closed $k$ of characteristic $p$.  But the countable
conjunction $\bigwedge_{K\in\N}\Psi_K$ just means

`` Given $|H|+1$ cycles in $Z_n(\tG(k))$, some two of them are
homologous. ''

So the cardinality of $H_n(\tG(k),R)$ is $|H|$ for all $k$.

Let now $k\ra K$ be an extension of algebraically closed fields.  The
induced map $H_n(\tG(k),R)\ra H_n(\tG(K),R)$ is injective (for example) by
model completeness of algebraically closed fields: if a cycle defined over
$k$ becomes a boundary over $K$, then a chain responsible for its being a
boundary must be definable already over $k$.  So within characteristic
$p$, all such maps must be isomorphisms, and
$H_n(\tG(k),R)=H_n(\tG(k_0),R)$ where $k_0$ is the algebraic closure of
the prime field.
\end{proof}

\begin{cor}
If $\tG$ satisfies an isoperimetric inequality in characteristic $p>0$,
then Friedlander's conjecture holds in that characteristic.
\end{cor}

This statement (which is the other half of Theorem A) follows by the
finiteness of $H_n(\tG(\overline{\FF}_p),\Z/l)$ and Lemma~\ref{2nd}. The
other half of Theorem B uses the first-order Lefschetz principle:

\begin{prop}  \label{milnor}
If the function \textup{(\asy)} of Theorem B exists for $\tG$, then Friedlander's conjecture holds.
\end{prop}

\begin{proof}
Fix $l$; for any $K\in\N$, set
$K_1=\max\{\asy(1),\asy(2),\dots,\asy(2K)\}$ and $|H|$ to be the common
cardinality of the groups $H_n(\tG(\overline{\FF}_p),\Z/l)$, $p\neq l$.
Consider the sentence $\Psi_K$ displayed above.  By assumption, it holds
over $k=\overline{\FF}_p$ for infinitely many primes $p$.  By the
Lefschetz principle, it holds over all algebraically closed fields of
characteristic zero.  That means that, in the notation of Lemma~\ref{1st},
the cardinality of $H_n(\tG(\CC)^\delta,\Z/l)$ is at most $|H_{\tG,n,l}|$.
Apply Lemma~\ref{1st}.
\end{proof}

A similar argument shows that the function $\asy$ must be at the same time
an isoperimetric bound for the group $\tG(\CC)^\delta$.  The converse
implication does not pass through the first-order Lefschetz principle,
and, unlike in the case of $\overline{\FF}_p$, no information is available
regarding $H_n(\tG(\overline{\QQ}),\Z/l)$ that would make
Prop.~\ref{oneway} applicable --- excepting those cases when the full
conjecture has already been proven!

\section{Stratifying the space of cycles}  \label{strat}
The following observation has long been known in saturated model theory,
but for convenience we include a proof.  By \emph{constructible subset} of
a variety $\V$ we mean one belonging to the boolean algebra generated by
Zariski-closed subsets of $\V$.  We only consider varieties defined over
some algebraically closed field $k$, and we identify them with their
$k$-points.

\begin{lemma}   \label{nice}
Let $\V$ be a variety over an uncountable, algebraically closed field $k$.
Suppose one has a collection $Z_i$, $i\in I$, of constructible subsets of
$\V$ such that $\card I<\card k$ and
\[
    Z := \bigcup_{i\in I} Z_i
\]
is constructible as well.  Then there exists a finite set
$i_1,i_2,\dots,i_N\in I$ such that
\[
   Z=Z_{i_1}\cup Z_{i_2}\cup\dots\cup Z_{i_N}.
\]
\end{lemma}

One proof of Lemma~\ref{nice} is akin to the `cylindrical' proof of the
compactness of the product of two compact topological spaces.  Without
loss of generality, we may assume $\V$ to be affine space $k^n$.  Also
without loss of generality, we may assume $Z=k^n$.  (Just add the
complement of $Z$ to the original collection.)

The proof is now by induction on $n$.  For $n=1$, the conclusion follows
since a constructible subset of $k$ is finite or co-finite, and by the
assumption $\card I<\card k$, one of the $Z_i$ have to be co-finite.
Assuming the claim holds below dimension $n$, write $n=r+s$ for some
$0<r,s<n$ and $k^n=A\times B$ with $A=k^r$, $B=k^s$.

For any $a\in A$, $\{a\}\times B$ is covered by its constructible subsets
$Z_i\cap(\{a\}\times B)$.  By the induction hypothesis for $k^s$, there
exists a finite index set $I_a\subseteq I$ such that
\[
   \{a\}\times B=\bigcup_{i\in I_a} Z_i\cap(\{a\}\times B).
\]
Each set defined as
\[
    C_a:=\{x\in A\;|\;\textup{\ for all\ }y\in B,\;
        \lc x,y\rc\in\bigcup_{i\in I_a} Z_i  \}
\]
forms a constructible subset of $A$, and their union is $A$.  As $a$
ranges over $A$, the range of $I_a$,
\[
   \{J\subseteq I\;|\;J=I_a\textup{\ for some\ }a\in A\}
\]
(thought of as a subset of the power set of $I$) has cardinality less than
that of $k$, since in fact the cardinality of all finite subsets of $I$
equals $\card I<\card k$.  By the induction hypothesis for $k^r$, one can
find finitely many $a_1,a_2,\dots,a_N\in A$ such that
$C_{a_1}\cup C_{a_2}\cup\dots\cup C_{a_N}=A$.  But that implies
\[
  k^n=A\times B=\bigcup_{\substack{j\in I_{a_i}\\i=1,2,\dots,N}} Z_j\,.
\]

\begin{prop}  \label{c1}
Let $k$ be an uncountable, algebraically closed field. If
$\card H_n(\tG(k),\Z/l) < \card k$, then $\tG(k)$ satisfies an
isoperimetric inequality.
\end{prop}

\begin{proof}
Having set up enough bookkeeping details, this becomes an immediate
consequence of Lemma~\ref{nice}.

\noindent
\textsl{Bookkeeping.}  From here on, $n$-chains will be thought of as
\emph{ordered} formal linear combinations of $n$-tuples of group elements.
Pick a representative $z_\alpha$ of each homology class $\alpha\in
H_n(\tG(k),\Z/l)$.  Let $0\in H_n(\tG(k),\Z/l)$ be represented by the
empty string.  Fix some size $K$ and tuple of coefficients $r_i\in\Z/l$,
$i=1,2,\dots,K$.  Define $\cZ_{r_1,r_2,\dots,r_K}$ as the locus in
$\tG(k)^{nK}$ of
$\big<\lc g_{11},g_{12},\dots,g_{1n}\rc,\dots,
\lc g_{K1},g_{K2},\dots,g_{Kn}\rc\big>$ such that
$\sum_{i=1}^K r_i\lc g_{i1},\dots,g_{in}\rc$ is a cycle in the bar
complex.  Define $Z(\alpha,K_1)$ as the locus in $\cZ_{r_1,r_2,\dots,r_K}$
of those cycles $z=\sum_{i=1}^K r_i\lc g_{i1},\dots,g_{in}\rc$ that
satisfy $z-z_\alpha=d_{n+1}(c)$ for some chain $c\in C_{n+1}(\tG(k))$ with
$\|c\|\leq K_1$.

$\cZ_{r_1,r_2,\dots,r_K}$ is an algebraic (i.e.\ Zariski-closed) subset of
$\tG(k)^{nK}$ and each $Z(\alpha,K_1)$ is a constructible subset of
$\cZ_{r_1,r_2,\dots,r_K}$.  Indeed, to say that
$z=\sum_{i=1}^K r_i\lc g_{i1},\dots,g_{in}\rc$ is a cycle is to say that
at least one of a finite number of possible cancellation patterns occurs
among the $(n+1)K$ many $n$-tuples that make up $d_{n+1}(z)$.  Each such
cancellation pattern is a system of equalities --- stated purely in terms
of the group multiplication on $\tG(k)$ --- among the group elements
$g_{ij}$.  Each $Z(\alpha,K_1)$ is a constructible subset of
$\cZ_{r_1,r_2,\dots,r_K}$.  Indeed, fix a $K_1$-tuple of coefficients for
the $n+1$-chain $c$, and denote by $[d_{n+1}(c)=z-z_\alpha]$ the locus in
$\tG(k)^{(n+1)K_1+nK}$ of those pairs $\lc c,z\rc$, $c\in
C_{n+1}(\tG(k))$, $z\in\cZ_{r_1,r_2,\dots,r_K}$ that satisfy
$d_{n+1}(c)=z-z_\alpha$.  $\big[d_{n+1}(c)=z-z_\alpha\big]$ is
Zariski-closed, and $Z(\alpha,K_1)$ is a finite union (as the $K_1$-tuple
of coefficients varies) of images of $\big[d_{n+1}(c)=z-z_\alpha\big]$
under the projection $\tG(k)^{(n+1)K_1+nK}\llra{\pr}\tG(k)^{nK}$.

Obviously $Z(\alpha,K_1)\subseteq Z(\alpha,K_2)$ for $K_1<K_2$, and
$Z(\alpha,K_1)\cap Z(\beta,K_2)=\nil$ for $\alpha\neq\beta$.  Since
\[
  \bigcup_{\substack{\alpha\in H_n(\tG(k),\Z/l)\\K_1\in\N}}
   Z(\alpha,K_1)=\cZ_{r_1,r_2,\dots,r_K}
\]
the lemma implies that one has a finite disjoint decomposition
\begin{equation}   \label{disj}
   \cZ_{r_1,r_2,\dots,r_K}=\bigsqcup_{j\in J} Z(\alpha_j,K_j).
\end{equation}
One (at most) of these homology classes $\alpha_j$, say $\alpha_0$, is the
zero one.  That means that every boundary of the form
$\sum_{i=1}^K r_i\lc g_{i1},\dots,g_{in}\rc$ possesses a filler of length
at most $K_0$.  Letting the tuple $r_1,r_2,\dots,r_K$ range over its
(finitely many!)\ possibilities, one obtains a finite value for
$\isop(K)$.
\end{proof}

Together with Cor.~\ref{same}, the next proposition completes the proof of
Theorem C of the introduction.

\begin{prop}   \label{c2}
If $\card H_n(\tG(k_0),\Z/l) < \card k_0$ for one uncountable,
algebraically closed field $k_0$, then, within the characteristic of
$k_0$, the groups $H_n(\tG(k),\Z/l)$ are countable for all algebraically
closed $k$, and every extension $k_1\ra k_2$ between algebraically closed
fields induces an isomorphism
$H_n(\tG(k_1),\Z/l)\llra{=}H_n(\tG(k_2),\Z/l)$.
\end{prop}

\begin{proof}
For brevity, let us introduce the notation $\dist(z_1,z_2)$ for
$\|z_1-z_2\|_\fil$ whenever $z_1,z_2$ are homologous cycles.  Note that
$\dist(z_1,z_3)\leq\dist(z_1,z_2)+\dist(z_2,z_3)$.  Over the field $k_0$,
one has a finite disjoint decomposition (\ref{disj}) of the space of
cycles $\cZ_{r_1,r_2,\dots,r_K}$ (corresponding to the tuple of
coefficients $r_i\in\Z/l$) into homology classes.  Without loss of
generality, we may assume that the cycle representatives $z_{\alpha_j}$
corresponding to the classes $\alpha_j$ occurring in the decomposition
themselves have size $K$.  Letting the tuple $r_i$ vary, one sees that for
each $K$, one can find $N$ and $K_1$ so that the sentence $\Phi_{K,N,K_1}$

`` There exist cycles $z_1,z_2,\dots,z_N$, all of size $K$, such that for
any cycle $z$ with $\|z\|=K$, either $\dist(z,z_1)\leq K_1$ or
$\dist(z,z_2)\leq K_1$ or $\dots$ or $\dist(z,z_N)\leq K_1$ ''

\noindent
holds over $k_0$.  But this is first-order, so by Tarski's theorem it
holds in all algebraically closed fields of the same characteristic as
$k_0$.

Note that the `obvious' thing to say (that the cycles $z_i$ are pairwise
non-homologous) is not first-order.  Nonetheless, for every $K$ there
exists a least $N=f(K)$ such that (for some $K_1<\infty$) $\Phi_{K,N,K_1}$
holds.  From the triangle inequality, one sees that over each
algebraically closed field, $f(K)$ is the number of distinct homology
classes that can be represented by cycles of size $K$; so this number does
not change with the underlying field.  A fortiori, $H_n(\tG(k),\Z/l)$ is
countable for every algebraically closed $k$ of the same characteristic as
$k_0$.

Let now $k_1\llra{i}k_2$ be a field extension as above, and consider the
induced $H_n(\tG(k_1),\Z/l)\llra{i_*}H_n(\tG(k_2),\Z/l)$.  It is always
injective, and if the sentences $\Phi_{K,N,K_1}$ hold, it is surjective
too.  Indeed, suppose $\alpha\in H_n(\tG(k_2),\Z/l)$ was not in the image
of $i_*$.  It would have to be represented by some $z\in Z_n(\tG(k_2))$,
say, of size $K$.  Since the inclusion $i$ does not change chain size,
this contradicts the injectivity of $i_*$ and the fact that the same
number of homology classes can be represented by cycles of size $K$ over
$k_1$ as over $k_2$.
\end{proof}

\section{The big picture}
Friedlander's conjecture concerns the effect of discretization $\CC^\delta\ra\CC$, and this paper revolves around the --- set-theoretically! --- equivalent
discretization $\prod_{\tP/\U}\overline{\FF}_p\ra\CC$.  It is natural to ask if (or in what sense) the two are compatible.  That amounts to pondering the
diagram
\begin{equation}\begin{split}  \label{bigpic}
\xymatrix{
H_n^\topo(B\tG(\CC)^\delta,\Z/l) \ar@{->>}[rr]^i &&
                                        H_n^\topo(B\tG(\CC),\Z/l) \\
H_n(\tG(\CC)^\delta,\Z/l) \ar@2{~}[r]^(.47)s\ar@{=}[u] &
\overset{\phantom{blah}}
{H_n(\prod_{\tP/\U}\tG(\overline{\FF}_p),\Z/l)}\ar@{->>}[r]^{[\iota]} &
\overset{\phantom{blah}}
{\prod_{\tP/\U}\,H_n(\tG(\overline{\FF}_p),\Z/l)} \ar@2{~}[u]_(.6)j }
\end{split}\end{equation}

Here $\tP$ is any infinite set of primes and $\U$ is any non-principal
ultrafilter on $\tP$.  $i$ is induced by the continuous homomorphism
$\tG(\CC)^\delta\ra\tG(\CC)$.  The left-hand vertical arrow is the
canonical isomorphism between discrete group homology and homology of
classifying spaces.  The isomorphism $s$ is induced by a \emph{choice} of
identification of $\CC^\delta$ with $\prod_{\tP/\U}\overline{\FF}_p$,
while $[\iota]$ is the canonical comparison map.  \emph{Choose} any
isomorphism $j_p$, independently for each $p\neq l$,
$H_*(\tG(\overline{\FF}_p),\Z/l)\approx H_*^\topo(B\tG(\CC),\Z/l)$; $j$ is
the ultraproduct of these isomorphisms, followed by the canonical
identification of $\prod_{\tP/\U}\,H_n^\topo(B\tG(\CC),\Z/l)$ with
$H_n^\topo(B\tG(\CC),\Z/l)$.

Friedlander proves that $i$ is surjective and conjectures that it is an isomorphism; Lemma~\ref{onto} proves that $[\iota]$ is surjective, and Theorem B states
that $[\iota]$ is an isomorphism if and only if $i$ is. It is \emph{highly} non-trivial, however, that the choices can be made compatibly so that this diagram
becomes commutative --- even up to isomorphism only.

Analyzing Lemma~\ref{onto}, one can show that for all $\tG$, $n$, and
(all but perhaps finitely many) $l$ there exist formulas
$z_1(-),z_2(-),\dots,z_{|H|}(-)$ in the language of rings such that for
any algebraically closed field $k$, the $z_1(k),z_2(k),\dots,z_{|H|}(k)$
are $n$-cycles in the bar complex of $\tG(k)$ with $\Z/l$ coefficients;
moreover, for almost all primes $p$, the cycles
$z_1(\overline{\FF}_p),z_2(\overline{\FF}_p),\dots,
z_{|H|}(\overline{\FF}_p)$ form exact representatives of
$H_n(\tG(\overline{\FF}_p),\Z/l)$.  As a corollary, one has, for almost
all $l$ and $p$, \emph{canonical} homomorphisms
\begin{equation}   \label{bigg}
  H_n(\tG(\overline{\FF}_p),\Z/l)\llra{h_p}
       H_n(\tG(\CC)^\delta,\Z/l)\llra{i}H_n^\topo(B\tG(\CC),\Z/l).
\end{equation}
Either of $h_p$ being an isomorphism or $i$ being an isomorphism is equivalent to Friedlander's conjecture; that the composite $j_p=i\circ h_p$ is an
isomorphism is Friedlander's theorem.  (Note that it is rather unobvious whether either Quillen's or Friedlander's proof of
$H_n(\tG(\overline{\FF}_p),\Z/l)\iso H_n^\topo(B\tG(\CC),\Z/l)$ gives a preferred isomorphism between the two sides.  The devil is in the passage between
positive and zero characteristics, which in Quillen's proof hinges on a Brauer lift, and in Friedlander's an embedding of the Witt vectors of
$\overline{\FF}_p$ in the complexes.)

Friedlander's conjecture implies that the inclusion $\overline{\QQ}\inc\CC$ induces an isomorphism
$H_n(\tG(\overline{\QQ}),\Z/l)\llra{=}H_n(\tG(\CC)^\delta,\Z/l)$, so the homology of $\tG(\CC)^\delta$ must have cycle representatives that are algebraic over
$\QQ$.  In fact, they must be $z_1(\overline{\QQ}),z_2(\overline{\QQ}),\dots,z_{|H|}(\overline{\QQ})$. The end result is that, under Friedlander's conjecture,
(\ref{bigpic}) becomes \emph{strictly} commutative for any choice of $\tP$, $\U$, and $s$ if the isomorphisms $j_p$ are chosen as in (\ref{bigg}).  Perhaps
this is the most beautiful embodiment of the compatibility of logic with geometry.


\providecommand{\bysame}{\leavevmode\hbox to3em{\hrulefill}\thinspace} \providecommand{\MR}{\relax\ifhmode\unskip\space\fi MR }
\providecommand{\MRhref}[2]{%
  \href{http://www.ams.org/mathscinet-getitem?mr=#1}{#2}
} \providecommand{\href}[2]{#2}

\end{document}